\documentclass{amsart}
\usepackage{color}
\usepackage{graphicx}
\usepackage{verbatim}
\usepackage{pstricks}
\newtheorem{thm}{Theorem}[section]

\newtheorem{lem}[thm]{Lemma}

\theoremstyle{definition}

\theoremstyle{remark}
\newtheorem{rem}[thm]{Remark}
\numberwithin{equation}{section}

\newcommand{\set}[1]{\left\{#1\right\}}
\newcommand{\Real}{\mathbb R}

\newcommand{\func}[1]{\ensuremath{\mathrm{#1} \:} }

\newcommand{\Area}[0]{\mathrm{Area}}
\newcommand{\Lip}[0]{\func{Lip}}

\title{Distortions of the Helicoid}
\author{Jacob Bernstein and Christine Breiner}
\begin{document}
\begin{abstract}
 Colding and Minicozzi have shown that an embedded minimal disk $0\in\Sigma\subset B_R$ in $\Real^3$ with large curvature at $0$ looks like a helicoid on the scale of $R$.  Near $0$, this can be sharpened: on the scale of $|A|^{-1}(0)$, $\Sigma$ is close, in a Lipschitz sense, to a piece of a helicoid.  We use surfaces constructed by Colding and Minicozzi to see this description cannot hold on the scale $R$.
 \end{abstract}
 \maketitle
\begin{figure}
 \includegraphics[width=3in]{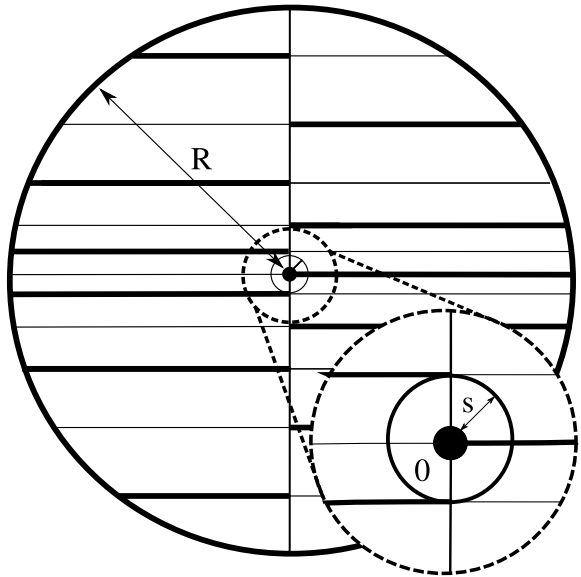}
 \caption{A cross section of one of Colding and Minicozzi's examples.  Here $R=1$ and $(0,s)$ is a blow-up pair.}\label{CMexample}
\end{figure}
\begin{figure}
 \includegraphics[width=3in]{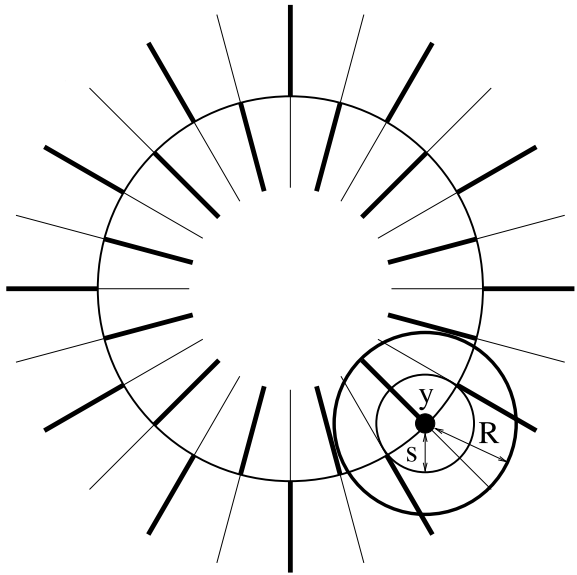}
 \caption{A cross section of one of Meeks and Weber's examples, with axis the circle.  Here $R$ is the outer scale of a disk and $(y,s)$ is a blow-up pair.}\label{MWexample}
\end{figure}

In \cite{CM1,CM2,CM3,CM4}, Colding and Minicozzi give a complete
description of the structure of embedded minimal disks in a ball in 
$\Real^3$.  Roughly speaking, they show that any such surface is
either modeled on a plane (i.e. is nearly graphical) or is modeled
on a helicoid (i.e. is two multi-valued graphs glued together
along an axis). In the latter case, the distortion may be quite
large.  For instance, in \cite{MW}, Meeks and Weber ``bend" the
helicoid; that is, they construct minimal surfaces where the axis
is an arbitrary $C^{1,1}$ curve (see Figure \ref{MWexample}).  A more serious example of
distortion is given by Colding and Minicozzi in \cite{CMPVNP}.
There they construct a sequence of minimal disks modeled on the
helicoid, but where the ratio between the scales (a measure of the
tightness of the spiraling of the multi-graphs) at different
points of the axis becomes arbitrarily large (see Figure \ref{CMexample}).  Note, locally, near
points of large curvature, the surface is close to a helicoid, and
so the distortions are necessarily global in nature.

Following \cite{CM2} we make the meaning of large curvature
precise by saying a pair $(y,s)\in\Sigma\times \Real^+$ is a $(C)$
\emph{blow-up pair} if $\sup_{B_s\cap \Sigma} |A|^2\leq
4C^2s^{-2}=4|A|^2(y)$ (here $C$ is large and fixed and $\Sigma
\subset \Real^3$ minimal). For $\Sigma$ minimal with $\partial
\Sigma\subset \partial B_R$ where $(0,s)$ is a blow-up pair, there
are two important scales; $R$ the outer scale and $s$ the blow-up
scale.  The work of Colding and Minicozzi gives a value
$0<\Omega<1 $ so that the component of $\Sigma \cap B_{\Omega R }$
containing $0$ consists of two multi-valued graphs glued together
(see for instance Lemma 2.5 of \cite{CY} for a self-contained
explanation).  On the other hand, Theorem 1.5 of \cite{BB} shows
that on the scale of $s$ (provided $R/s$ is large), $\Sigma$ is
bi-Lipschitz to a piece of a helicoid with Lipschitz constant near
1.  Using the surfaces constructed in \cite{CMPVNP} we show that
such a result cannot hold on the outer scale and indeed fails to hold on certain smaller scales:
\begin{thm} \label{FirstThm}
Given $1>\Omega,\epsilon>0$ and $1/2 > \gamma \geq 0$ there exists an embedded
minimal disk $0\in\Sigma$ with $\partial \Sigma\subset
\partial B_{R}$ and $(0,s)$ a blow-up pair so:  the component of
$B_{\Omega R^{1-\gamma} s^{\gamma}}\cap \Sigma$ containing $0$ is not
bi-Lipschitz to a piece of a helicoid with Lipschitz constant in
$((1+\epsilon)^{-1},1+\epsilon)$.
\end{thm}
First, we recall the surfaces constructed in \cite{CMPVNP}:
\begin{thm} \label{SurfExst}(Theorem 1 of \cite{CMPVNP}) There is a sequence of compact embedded minimal disks $0\in \Sigma_i\subset B_1\subset \Real^3$ with $\partial \Sigma_i\subset \partial B_1$ containing the vertical segment $\set{(0,0,t) : |t|\leq 1}\subset \Sigma_i$ such that the following conditions are satisfied:
\begin{enumerate}
\item $\lim_{i\to \infty} |A_{\Sigma_i}|^2(0)\to \infty$
\item \label{SecItema}
$\sup_{\Sigma_i}|A_{\Sigma_i}|^2 \leq 4|A_{\Sigma_i}|^2(0) = 8
a_i^{-4}$ for a sequence $a_i \to 0$ 
\item
\label{SecItem} $\sup_i \sup_{\Sigma_\backslash B_\delta}
|A_{\Sigma_i}|^2 <K\delta^{-4}$ for all $1>\delta>0$ and $K$ a
universal constant. 
\item \label{ThrItem}
$\Sigma_i\backslash\set{x_3-axis} =\Sigma_{1,i}\cup \Sigma_{2,i}$
for multi-valued graphs $\Sigma_{1,i}$ and $\Sigma_{2,i}.$
\end{enumerate}
\end{thm}
\begin{rem}
\eqref{SecItema} and \eqref{SecItem} are slightly sharper than what is stated in Theorem
1 of \cite{CMPVNP}, but follow easily. \eqref{SecItema} follows from the Weierstrass data (see Equation (2.3) of \cite{CMPVNP}). This also gives \eqref{SecItem} near the axis, whereas away from the
axis use \eqref{ThrItem} and Heinz's curvature estimates.
\end{rem}
Next introduce some notation.  For a surface $\Sigma$ (with a smooth metric)  we
denote intrinsic balls by $\mathcal{B}_s^\Sigma$ and define the
(\emph{intrinsic}) \emph{density ratio} at a point $p$ as: $
\theta_s (p,\Sigma)=(\pi s^2)^{-1}\Area(\mathcal{B}_s^\Sigma (p)
)$. When $\Sigma$ is immersed in $\Real^3$ and has the induced metric,
$\theta_s(p,\Sigma)\leq \Theta_s (p,\Sigma)=({\pi
s^2})^{-1}{\Area(B_s(p)\cap \Sigma)}$, the usual (extrinsic)
density ratio.  Importantly, the intrinsic density ratio is well-behaved under bi-Lipschitz maps.  Indeed, if $f: \Sigma \to
\Sigma'$ is injective and with $\alpha^{-1}<\Lip f<\alpha$, then:
\begin{equation} \label{DensBnds}
\alpha^{-4}\theta_{\alpha^{-1} s} (p,\Sigma) \leq \theta_s (f(p),{\Sigma'}) \leq \alpha^4\theta_{\alpha s} (p,\Sigma).
\end{equation}
This follows from the inclusion,
$\mathcal{B}^\Sigma_{\alpha^{-1} s} (f^{-1}(p))\subset f^{-1}
(\mathcal{B}^{\Sigma'}_s (p))$ and the behavior of area under
Lipschitz maps, $\Area( f^{-1} (\mathcal{B}^{\Sigma'}_s (p))\leq
(\Lip f^{-1})^2 \Area (\mathcal{B}^{\Sigma'}_s (p))$.

Note that by standard area estimates for minimal graphs, if $\Sigma\cap B_s(p)$ is a minimal graph then
$\theta_s(p,\Sigma)\leq 2$.  In contrast, for a point near the
axis of a helicoid, for large $s$ the density ratio is large.
Thus, in a helicoid the density ratio for a fixed, large $s$
measures, in a rough sense, the distance to the axis.  More generally, this holds near blow-up pairs of embedded minimal disks:
\begin{lem} \label{LowerBndDensLem}
Given $D>0$ there exists $R>1$ so: If $0\in\Sigma\subset B_{2Rs}$
is an embedded minimal disk with $\partial \Sigma\subset\partial
B_{2 Rs}$ and $(0,s)$ a blow-up pair then
$\theta_{Rs}(0,\Sigma)\geq D$.
\end{lem}
\begin{proof}
We proceed by contradiction, that is suppose there were a $D>0$
and embedded minimal disks $0\in\Sigma_i$ with
$\partial \Sigma_i \subset \partial
B_{2 R_i s}$ with $R_i\to \infty$ and $(0,s)$ a blow-up pair so
that $\theta_{R_i s}(0,\Sigma_i)\leq D$.  The chord-arc
bounds of \cite{CY} imply there is a $1>\gamma>0$ so $\mathcal{B}_{R_i
s}^{\Sigma_i} (0)\supset \Sigma_i \cap B_{\gamma R_i s}$.  Hence,
the
intrinsic density ratio bounds the extrinsic density ratio, i.e. $D\geq \theta_{R_i s}(p,{\Sigma_i})\geq
\gamma^2 \Theta_{\gamma R_i s}(p, \Sigma_i)$. Then, by a result of
Schoen and Simon \cite{ScSi} there is a constant
$K=K(D\gamma^{-2})$, so $|A_{\Sigma_i}|^2(0)\leq K(\gamma R_i
s)^{-2}$.  But for $R_i$ very large this contradicts that $(0,s)$
is a blow-up pair for all $\Sigma_i$.
\end{proof}
\begin{rem}
Note that the above does not depend on the strength of chord-arc
bounds. In fact, it is also an immediate consequence of the fact
that intrinsic area bounds on a disk give total curvature bounds.
In turn, the total curvature bounds again yield uniform curvature
bounds. See Section 1 of \cite{CM2} for more detail.
\end{rem}

To produce our counterexample, we exploit the fact that two points
on a helicoid that are equally far from the axis must have the
same density ratio.  Assuming the existence of a Lipschitz map
between our surface $\Sigma$ and a helicoid, we get a
contradiction by comparing the densities for two appropriately
chosen points that map to points equally far from the axis of the
helicoid.
\begin{proof}(of Theorem \ref{FirstThm})
Fix $1>\Omega,\epsilon>0$ and $1/2 > \gamma\geq 0$ and set
$\alpha=1+\epsilon$. Let $\Sigma_i$ be the surfaces of Theorem
\ref{SurfExst}; we claim for $i$ large, $\Sigma_i$ will be the
desired example. Suppose this was not the case. Setting $s_i =
Ca_i^2/\sqrt{2}$, where $a_i$ is as in \eqref{SecItema} and $C$ is the blow-up constant,one has
$(0,s_i)$ is a blow-up pair in $\Sigma_i$, since $\sup_{\Sigma_i
\cap B_{s_i}} |A_{\Sigma_i}|^2\leq 8a_i^{-4} = 4C^2s_i^{-2} =
4|A_{\Sigma_i}|^2(0)$, moreover, $s_i\to 0$. Hence, with $R_i=
\Omega s_i^{\gamma}<1$, the component of $B_{R_i}\cap {\Sigma}_i$
containing $0$, ${\Sigma}_i'$, is bi-Lipschitz to a piece of a
helicoid with Lipschitz constant in $(\alpha^{-1},\alpha)$. That
is, there are subsets $\Gamma_i$ of helicoids and diffeomorphisms
$f_i:{\Sigma}_i' \to \Gamma_i$ with $\Lip f_i\in
(\alpha^{-1},\alpha)$.

We now begin the density comparison.  First, Lemma
\ref{LowerBndDensLem} implies there is a constant $r>0$ so for $i$
large $\theta_{r s_i} (0,{\Sigma_i'})\geq 4 \alpha^8$ and thus by
\eqref{DensBnds} $\theta_{\alpha r s_i} (f_i(0),{\Gamma_i})\geq 4
\alpha^{4} $. We proceed to find a point with small density on
$\Sigma_i$ that maps to a point on $\Gamma_i$ equally far from the
axis as $f_i(0)$ (which has large density).

Let $U_i$ be the (interior) of the component of $B_{1/2 R_i} \cap
{\Sigma}_i$ containing $0$. Note for $i$ large enough, as $s_i/R_i\to 0$, the distance between $\partial U_i$ and $\partial {\Sigma_i}'$
is greater than $4\alpha^2 rs_i$. Similarly, for $p\in\partial U_i$ for $i$ large, $p'\in
\mathcal{B}^{\Sigma_i'}_{4\alpha^2 r s_i}(p)$ implies $|p'|\geq
\frac{1}{4} R_i$.  Hence, property \eqref{SecItem} gives that
$|A_{\Sigma_i'}|^2(p')\leq K' s_i^{-4\gamma}$.  Thus, for $i$
sufficiently large $\mathcal{B}_{\alpha^2 r s_i} (p)$ is a graph
and so $\theta_{\alpha^2 r s_i} (p,{{\Sigma}_i'})\leq 2$. Pick
$u_i\in
\partial f(U_i)$ at the same distance to the axis as $f_i(0)$ and so the density ratio is the same at both points (see Figure \ref{ProofImg}).
As $f_i(U_i)$ is an open
subset of $\Gamma_i$ containing $f_i(0)$, ${p}_i=f_i^{-1}(u_i)\in
\partial U_i$. Notice that $\theta_{\alpha  r s_i}
(u_i,{\Gamma_i})=\theta_{\alpha r s_i} (f_i(0),{\Gamma_i})\geq 4
\alpha^{4}$ so $2 \alpha^4 \geq\alpha^4 \theta_{\alpha^2 r
s_i}(p_i,{{\Sigma}_i'})\geq  4 \alpha^{4}$.
\begin{figure}
 \includegraphics[width=5in]{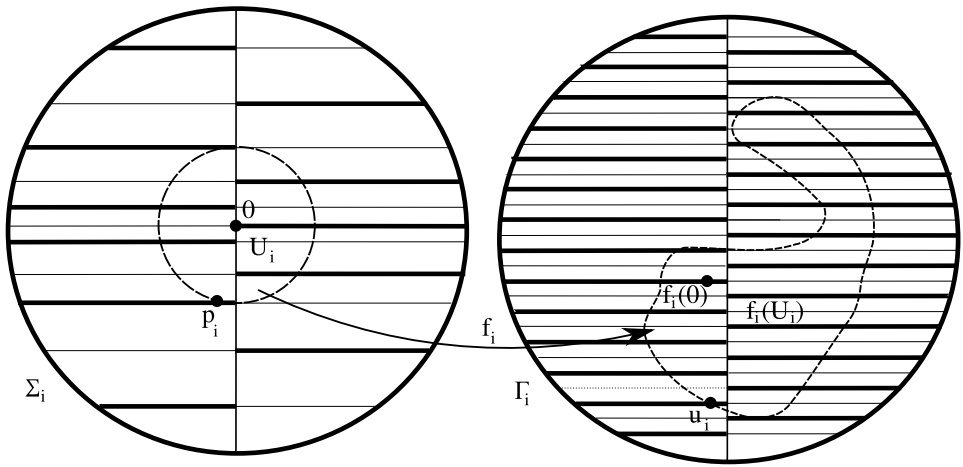}
 \caption{Finding $u_i$ }\label{ProofImg}
\end{figure}

\end{proof}

\bibliographystyle{amsplain}
\bibliography{FinalDraft}

\providecommand{\bysame}{\leavevmode\hbox to3em{\hrulefill}\thinspace}
\providecommand{\MR}{\relax\ifhmode\unskip\space\fi MR }
\providecommand{\MRhref}[2]{%
  \href{http://www.ams.org/mathscinet-getitem?mr=#1}{#2}
}
\providecommand{\href}[2]{#2}
\begin{thebibliography}{1}

\bibitem{BB}
J.~Bernstein and C.~Breiner, \emph{Helicoid-like minimal disks and uniqueness},
  Preprint.

\bibitem{CMPVNP}
T.~H. Colding and W.~P.~Minicozzi II, \emph{{Embedded minimal disks: Proper
  versus nonproper-Global versus local}}, Trans. AMS \textbf{356} (2003),
  no.~1, 283--289.

\bibitem{CM1}
\bysame, \emph{{The space of embedded minimal surfaces of fixed genus in a
  3-manifold I; Estimates off the axis for disks}}, Ann. of Math. (2)
  \textbf{160} (2004), no.~1, 27--68.

\bibitem{CM2}
\bysame, \emph{{The space of embedded minimal surfaces of fixed genus in a
  3-manifold II; Multi-valued graphs in disks}}, Ann. of Math. (2) \textbf{160}
  (2004), no.~1, 69--92.

\bibitem{CM3}
\bysame, \emph{{The space of embedded minimal surfaces of fixed genus in a
  3-manifold III; Planar domains}}, Ann. of Math. (2) \textbf{160} (2004),
  no.~2, 523--572.

\bibitem{CM4}
\bysame, \emph{{The space of embedded minimal surfaces of fixed genus in a
  3-manifold IV; Locally simply connected}}, Ann. of Math. (2) \textbf{160}
  (2004), no.~2, 573--615.

\bibitem{CY}
\bysame, \emph{{The Calabi-Yau conjectures for embedded surfaces}}, Ann. of
  Math. \textbf{167} (2008), no.~1, 211--243.

\bibitem{MW}
W.H. Meeks and M.~Weber, \emph{{Bending the helicoid}}, Mathematische Annalen
  \textbf{339} (2007), no.~4, 783--798.

\bibitem{ScSi}
R.~Schoen and L.~Simon, \emph{{Regularity of simply connected surfaces with
  quasiconformal Gauss map}}, Seminar on minimal submanifolds, Ed. E. Bombieri,
  Ann. of Math. Studies \textbf{103}, 127--145.

\end{thebibliography}
\end{document}